\documentclass[journal]{IEEEtran}
\def\baselinestretch{1.2}

\usepackage{tikz}
\graphicspath{{images/}}
\usepackage{amsmath}
\usepackage{latexsym, amssymb}
\usepackage{graphicx}
\usepackage{amsmath, amsbsy}
\usepackage{amsopn, amstext}
\usepackage{ifpdf,hyperref}
\usepackage{cancel, color}
\usepackage{epstopdf}

\def\lplus{\,\rotatebox[]{-90}{$\pm$}\,}
\def\rplus{\,\rotatebox[]{90}{$\pm$}\,}
\def\lminus{\vdash}
\def\rminus{\dashv}

\def\GL{GL}
\def\gl{gl}

\def\End{End}
\def\ad{ad}

\DeclareMathOperator{\trace}{trace}

\def\cal{\mathcal}

\def\Tr{Tr}
\def\trace{trace}

\def\diag{diag}
\def\ra{\rightarrow}

\def\a{\alpha}
\def\b{\beta}
\def\d{\delta}

\def\0{{\bf 0}}

\def\A{\left<A\right>}
\def\B{\left<B\right>}
\newcommand{\R}{{\mathbb R}}
\newcommand{\Q}{{\mathbb Q}}
\newcommand{\C}{{\mathbb C}}

\newtheorem{thm}{Theorem}[section]
\newtheorem{dfn}[thm]{Definition}
\newtheorem{prp}[thm]{Proposition}
\newtheorem{exa}[thm]{Example}
\newtheorem{lem}[thm]{Lemma}

\newtheorem{rem}[thm]{Remark}

\begin{document}

\title{Dimension-Free Matrix Spaces}

\author{Daizhan Cheng  %,Zhengping Ji
	\thanks{This work is supported partly by the National Natural Science Foundation of China (NSFC) under Grants 62073315, 61074114, and 61273013.}
	\thanks{Key Laboratory of Systems and Control, Academy of Mathematics and Systems Sciences, Chinese Academy of Sciences,
		Beijing 100190, P. R. China (e-mail: dcheng@iss.ac.cn).}%, jizhengping@amss.ac.cn).}
    %\thanks{Zequn Liu and Hongsheng Qi are with the Key Laboratory of Systems and Control, Academy of Mathematics and Systems Sciences, Chinese Academy of Sciences, Beijing 100190, and University of Chinese Academy of Sciences, Beijing 100049, P. R. China (e-mail: liuzequn@amss.ac.cn, qihongsh@amss.ac.cn).}
%    \thanks{Corresponding author: Daizhan Cheng. Tel.: +86 10 82541232.}
}

%\markboth{IEEE Transactions On Automatic Control, Vol. XX, No. Y, Month 202Z}
%{Cheng \MakeLowercase{\textit{et al.}}: Field Extension via Semi-tensor Product of Matrices}

\maketitle

\begin{abstract}
Based on various types of semi-tensor products of matrices, the corresponding equivalences of matrices are proposed. Then the corresponding vector space structures are obtained as the quotient spaces under equivalences, which are called the dimension-free Matrix spaces (DFESs). Certain structures and properties are investigated. Finaly, the Lie bracket structure of general linear algebra is extended to DfMSs to make them Lie algebras, called dimension-free general linear algebra (DFGLA).  Inspire of the fact that the DFGLAs are of  infinite dimension, they have most properties of finite dimensional Lie algebras, whicl are studied in the paper.
\end{abstract}


\IEEEpeerreviewmaketitle

\section{Introduction}

Lie group was firstly proposed by Sophus Lie (1840-1899) to describe the continuous symmetries of algebraic or geometric objects, and the corresponding Lie algebra is considered as the infinitesimal generator of the Lie group. Since the objects might be of finite or infinite dimensions, the related Lie groups and Lie algebras can also be of either finite or infinite dimensions. As pointed in \cite{kac83}: ``Starting with the works of Lie, Killing, and Cartann, the theory of finite-dimensional Lie groups and Lie algebras has developed systematically in depth and scope. On the other hand, works on simple infinite dimensional Lie algebras have been virtually forgotten until the mod-sixties" (of last century). By the author's observation, the reason for this may lies on the big differences between Lie groups and Lie algebras of finite vs infinite dimensions.

Finite dimensional Lie algebra has many nice properties, such as Killing form, complete reducibility of representations, root systems, etc. \cite{hum72,wan13}, which can not be easily extended to infinite dimensional case.

The STP of matrices was proposed firstly by the first author as follows:
Let $A\in {\cal M}_{m\times n}$ and $B\in {\cal M}_{p\times q}$, $t=n\vee p$ be the least common multiple of $n$ and $p$. The STP of $A$ and $B$ is defined as follows:
\begin{align}\label{1.1}
A\ltimes B:=\left(A\otimes I_{t/n}\right)\left(B\otimes I_{t/p}\right),
\end{align}
where $\otimes$ is the Kronecker product.

STP generalizes the conventional matrix product to general case where two factor matrices are of arbitrary dimensions \cite{che11,che12}.  Fortunately, STP keeps major properties of the conventional matrix product available. The last two decades have witnessed a fast progress in the theory and applications of STP of matrices. We refer to some survey papers for its some major applications: (i) application to logical dynamic (control) networks \cite{for16,lu17,muh16}; (ii) application to game theory \cite{che21}; (iii) application to finite automata \cite{yanpr}; (iv) application to engineering problems \cite{li18}.

Based on STP, certain equivalences of matrices are obtained. Then the quotient spaces are produced as infinite dimensional vector space of matrices, called DFMSs. The addition and scalar product are posed on quotient space to make it a vector space.
The lattice and the topology are also obtained. The structure of general linear algebras (GLAs) ($\gl(n,\R)$ or $\gl(n,\C)$) can be extended to DFMSs, which turns the DFMS into an infinite dimensional Lie algebra, called the DFGLA. Then we prove that many properties of GLA can be transferred to DFGLA. Finally, the Killing structure is also built for DFMSs.

The outline of this paper is as follows: In section 2 several general STPs of matrices are presented by using matrix multiplier. Section 3 poses addition and scalar product of matrices to make certain sets of matrices   vector spaces. In Section 4 a kind of equivalence relation is proposed for matrices, which lead to the lattice structure on the set of matrices. The quotient spaces, called the DFMSs, are obtained by the equivalences, which are discussed in Section 5. Section 6 poses the Lie algebraic structure onto the DFMSs to build an infinite dimensional Lie algebra on DFMS. The Killing form is introduced for the DFMSs in Section 7. Section 8 is a brief conclusion.

\section{General STP of Matrices}

The STP of matrices defined by (\ref{1.1}) is the classical one. There are some other useful STPs \cite{che19}. A general way to construct STPs is via matrix multiplier, which is defined as follow:

\begin{dfn}\label{d2.1}
A set of nonzero matrices
$$
\varGamma:=\left\{\varGamma_n\in {\cal M}_{n\times n}\;|\; n\geq 1\right\}
$$
is called a matrix multiplier, if it satisfies:
\begin{itemize}
\item[(i)]~~
\begin{align}\label{2.1}
\varGamma_n\varGamma_n=\varGamma_n;
\end{align}
\item[(ii)]~~
\begin{align}\label{2.2}
\varGamma_p\otimes \varGamma_q=\varGamma_{pq}.
\end{align}
\end{itemize}
\end{dfn}

It follows from the definition immediately that the matrix multiplier has some elementary properties as follows.

\begin{prp}\label{p2.2} Consider a matrix multiplier  $\varGamma$.
\begin{itemize}
\item[(i)] The spectrum (set of eigenvalues) of ~$\varGamma_n$ is either ~$\sigma(\varGamma_n)=\{1\}$ or ~$\sigma(\varGamma_n)=\{0,1\}$. Hence
\begin{align}\label{2.3}
\sigma_{\max}(\varGamma_n)=1,\quad n=1,2,\cdots
\end{align}
\item[(ii)]
\begin{align}\label{2.4}
\varGamma_1=1.
\end{align}
\end{itemize}
\end{prp}

Using a matrix multiplier, the corresponding STP of matrices can be defined as following.

\begin{dfn}\label{d2.3} Given a matrix multiplier $\varGamma$, and let $A\in {\cal M}_{m\times n}$ and $B\in {\cal M}_{p\times q}$, $t=n\vee p$ be the least common multiple of $n$ and $p$.
\begin{itemize}
\item[(i)] The left $\varGamma$-STP of $A$ and $B$ is defined by
\begin{align}\label{2.5}
A\ltimes_{\varGamma} B:=\left(A\otimes \varGamma_{t/n}\right)\left(B\otimes \varGamma_{t/p}\right).
\end{align}
\item[(ii)] The right $\varGamma$-STP of $A$ and $B$ is defined by
\begin{align}\label{2.6}
A\rtimes_{\varGamma} B:=\left(\varGamma_{t/n}\otimes A\right)\left(\varGamma_{t/p}\otimes B\right).
\end{align}
\end{itemize}
\end{dfn}

Using the properties of $\varGamma$, it is easy to verify the following proposition.

\begin{prp}\label{p2.4} Both left and right $\varGamma$-STP, defined by (\ref{2.5}) and (\ref{2.6}) respectively, satisfy the following condition.
\begin{itemize}
\item[(i)] If $n=p$, then
$$
A\ltimes_{\varGamma} B=A\rtimes_{\varGamma} B=AB.
$$
That is, they are the generalizations of the convenient matrix product.

\item[(ii)] (Associativity)

\begin{align}\label{2.7}
(A\bowtie_{\varGamma} B)\bowtie_{\varGamma} C =A \bowtie_{\varGamma} (B \bowtie_{\varGamma} C),
\end{align}
where $\bowtie_{\varGamma}$ can be either $\ltimes_{\varGamma}$ or $\rtimes_{\varGamma}$.

\item[(iii)] (Distributivity)

\begin{align}\label{2.8}
\begin{array}{l}
(A+B) \bowtie_{\varGamma} C:=A\bowtie_{\varGamma} C+B\bowtie_{\varGamma} C\\
A\bowtie_{\varGamma} (B+C)=A\bowtie_{\varGamma} B+A\bowtie_{\varGamma} C.
\end{array}
\end{align}

\item[(iv)] (Transpose)

\begin{align}\label{2.9}
(A\bowtie_{\varGamma} B)^T=B^T\bowtie_{\varGamma} A^T.
\end{align}

\item[(v)] (Inverse)

Assume ~$\varGamma_n$, ~$n\geq 1$ are invertible. If ~$A$ and ~$B$ are invertible, then ~$A\bowtie_{\varGamma} B$ is invertible. Moreover,
\begin{align}\label{2.10}
(A\bowtie_{\varGamma} B)^{-1}=B^{-1}\bowtie_{\varGamma} A^{-1}.
\end{align}
\end{itemize}
\end{prp}

The following example shows some matrix multipliers.

\begin{exa}\label{e2.5}
\begin{itemize}
\item[(i)]  (MM-1 STP): Set
$$
\varGamma^1:=\{I_n\;|\;n=1,2,\cdots\}.
$$
Then we denote
$$
\ltimes_{\varGamma}=\ltimes;\quad \rtimes_{\varGamma}=\rtimes.
$$
This set is called the type-1 STP. Eq. (\ref{1.1}) is of this type.

\item[(ii)] (MM-2 STP): Set
$$
\varGamma^2:=\{J_n\;|\;n=1,2,\cdots\},
$$
where
\begin{align}\label{2.11}
J_n:=\frac{1}{n}{\bf 1}_{n\times n},\quad n=1,2,\cdots.
\end{align}
Then we denote
$$
\ltimes_{\varGamma}=\circ_{\ell};\quad \rtimes_{\varGamma}=\circ_r.
$$
This set is called the type-2 STP.

\item[(iii)]  Set
$$
\varGamma^3:=\{U_n\;|\;n=1,2,\cdots\},
$$
where,
\begin{align}\label{2.12}
U_n:=\begin{bmatrix}
1&0&\cdots&0\\
0&0&\cdots&0\\
\vdots&~&~&~\\
0&0&\cdots&0\\
\end{bmatrix}\in {\cal M}_{n\times n}.
\end{align}

\item[(iv)]  Set
$$
\varGamma^4:=\{L_n\;|\;n=1,2,\cdots\},
$$
where
\begin{align}\label{2.13}
L_n:=\begin{bmatrix}
0&0&\cdots&0\\
0&0&\cdots&0\\
\vdots&~&~&~\\
0&0&\cdots&1\\
\end{bmatrix}\in {\cal M}_{n\times n}.
\end{align}
\end{itemize}
\end{exa}

\section{Vector Space Structure on Matrices}

Denote by
$$
{\cal M}:=\bigcup_{m=1}^{\infty}\bigcup_{n=1}^{\infty} {\cal M}_{m\times n}.
$$
Then the following proposition is obvious.

\begin{prp}\label{p3.1} Consider the product $\bowtie_{\varGamma}$ on ${\cal M}$. Then, $\left({\cal M},\bowtie_{\varGamma}\right)$ is a semi-group. Particularly, set $\varGamma=\varGamma^1=\{I_n\;|\;n=1,2,\cdots\}$, then  $\left({\cal M},\bowtie_{\varGamma}\right)$ is a monoid, i.e., a semi-group with identity.
\end{prp}

Define
$$
{\cal M}_{\mu}:=\{A\in {\cal M}_{m\times n}\;|\;m/n=\mu\}.
$$
Then we have the following partition:
\begin{align}\label{3.1}
{\cal M}=\bigcup_{\mu\in \Q_+}{\cal M}_{\mu},
\end{align}
where $\Q_+$ is the set of positive natural numbers.

Assume
$$
\mu=\mu_y/\mu_x,
$$
where $\mu_y\wedge \mu_x=1$, i.e., $\mu_y$ and $\mu_x$ are co-prime. Then $\mu_y$ and $\mu_x$ are called the $y$ and $x$ component of $\mu$ respectively. Hereafter, we always assume $\mu_y$ and $\mu_x$ are co-prime. Then they are uniquely determined by $\mu$.

Next, we pose a vector space structure on ${\cal M}_{\mu}$. To this end, we define the addition on it.

\begin{dfn}\label{d3.2} Consider ${\cal M}_{\mu}$, where $\mu\in \Q_+$. Let $A,B\in {\cal M}_{\mu}$ and
$A\in {\cal M}_{m\times n}$,  $B\in {\cal M}_{p\times q}$ respectively, and $t=m\vee p$.
Then
\begin{itemize}
\item[(i)] the left $\varGamma$ addition is defined by
\begin{align}\label{3.2}
A\lplus_{\varGamma}B:=\left(A\otimes \varGamma_{t/m}\right)
+\left(B\otimes \varGamma_{t/p}\right);
\end{align}
\item[(ii)] the right $\varGamma$ addition is defined by
\begin{align}\label{3.3}
A\rplus_{\varGamma}B:=\left(\varGamma_{t/m}\otimes A\right)
+\left(\varGamma_{t/p}\otimes B\right).
\end{align}
\end{itemize}
\end{dfn}

Corresponding to (\ref{3.2}) and (\ref{3.3}), the subtraction can be defined respectively as
\begin{align}\label{3.4}
A\lminus_{\varGamma}B:=A\lplus_{\varGamma}(-B);
\end{align}
and
\begin{align}\label{3.5}
A\rminus_{\varGamma}B:=A\rplus_{\varGamma}(-B).
\end{align}

By definition, it is easy to verify the following.

\begin{prp}\label{p3.3} With the addition defined by (\ref{3.2}) (or (\ref{3.3})) and a conventional scalar product, ${\cal M}_{\mu}$ is a pseudo vector space.
\end{prp}

\begin{rem}\label{r3.4} A pseudo vector space satisfies all the requirements for a vector space except one:
the ``zero" is not unique \cite{abr78}. Hence the inverse $-x$ for each $x$ is also not unique. In our case,
$$
\vec{\bf 0}:=\{{\bf 0}_{k\mu_y\times k\mu_x}\;|\;k=1,2,\cdots\}.
$$
\end{rem}

\section{Equivalence and Lattice Structure on ${\cal M}_{\mu}$}

\begin{dfn}\label{d4.1} Let $A,B\in {\cal M}_{\mu}$.
\begin{itemize}
\item[(i)]
$A$ is said to be left $\varGamma$- equivalent to $B$, denoted by $A\sim^{\varGamma}_{\ell} B$, if there exist $\varGamma_{\a}$ and $\varGamma_{b}$ such that
\begin{align}\label{4.1}
A\otimes \varGamma_{\a}=B\otimes \varGamma_{\b}.
\end{align}
\item[(ii)]
$A$ is said to be right $\varGamma$- equivalent to $B$, denoted by $A\sim^{\varGamma}_{r} B$, if there exist $\varGamma_{\a}$ and $\varGamma_{b}$ such that
\begin{align}\label{4.2}
\varGamma_{\a}\otimes A=\varGamma_{\b}\otimes B.
\end{align}
\end{itemize}
\end{dfn}

\begin{rem}\label{r4.2} Hereafter, we consider only
\begin{itemize}
\item Case 1:
$$
\varGamma=\varGamma^1=\{I_n\;|\;n=1,2,\cdots\}.
$$
Then $A\sim^{\varGamma^1}_{\ell}B$ is simply denoted by
$$
A\sim_{\ell}B.
$$
The equivalence class is denoted by
$$
\left<A\right>_{\ell}:=\{B\;B\sim_{\ell}A\}.
$$

$A\sim^{\varGamma^1}_{r}B$ is simply denoted by
$$
A\sim_{r}B.
$$
The equivalence class is denoted by
$$
\left<A\right>_{r}:=\{B\;|\; B\sim_{r}A\}.
$$

\item Case 2:
$$
\varGamma=\varGamma^2=\{J_n\;|\;n=1,2,\cdots\}.
$$
Then $A\sim^{\varGamma^2}_{\ell}B$ is simply denoted by
$$
A\approx_{\ell}B.
$$
The equivalence class is denoted by
$$
\ll A\gg_{\ell}:=\{B\;|\;B\approx_{\ell} A\}.
$$

$A\sim^{\varGamma^2}_{r}B$ is simply denoted by
$$
A\approx_{r}B.
$$
The equivalence class is denoted by
$$
\ll A\gg_{r}:=\{B\;|\;B\approx_{r}A\}.
$$
\end{itemize}
\end{rem}

For statement ease, in the following we consider only $\varGamma=\varGamma^1$ and the equivalence
$$
\sim:=\sim^{\varGamma^1}_{l}.
$$
Hence
$$
\left<A\right>:=\left<A\right>_{\ell}.
$$
All the arguments are applicable to the other three cases.

\begin{dfn}\label{d4.3}  Consider $\left<A\right>$. Assume $A,B\in \left<A\right>$, if there exists an $I_{\a}$ such that $A\otimes I_{\a}=B$, then $A$ is said to be preceding to $B$, denoted by $A\prec B$.
Then $\left<A\right>$ becomes a partially ordered set.
\end{dfn}

\begin{prp}\label{p4.4} \cite{che19b} Let $A,B\in {\cal M}_{\mu}$. Assume $A\sim B$, that is, (\ref{4.1}) holds. Then there exist a $\Lambda\in {\cal M}_{\mu}$ and two identity matrices $I_{\a}$ and $I_{\b}$, such that
\begin{align}\label{4.3}
A=\Lambda\otimes I_{\b},\quad B=\Lambda\otimes I_{\a}.
\end{align}
\end{prp}

Consider (\ref{4.1}). Without loss of generality, we can assume $\a\wedge \b=1$. In this case we define
\begin{align}\label{4.4}
\Theta:=A\otimes I_{\a}=B\otimes I_{\b}.
\end{align}

\begin{prp}\label{p4.5} \cite{che19b}  Assume $\a \wedge \b=1$. Then the least common upper boundary of $A$ and $B$ is $\Theta$, and the greatest common lower boundary of $A$ and $B$ is $\Lambda$.
\end{prp}

\begin{dfn}\label{d4.6} \cite{bur81} Consider a non-empty set $L\neq \emptyset$, if there is a partial order $\prec$ on it such tha for any two elements $a,b$ in $L$ there exist a least common upper boundary $\inf(a,b)\in L$ and a greatest common
lower boundary $\sup(a,b)\in L$, then $(L,\prec)$ is called a lattice.
\end{dfn}

Recall Proposition \ref{p4.5}, one sees easily that $(\left<A\right>, \prec)$ is a lattice.

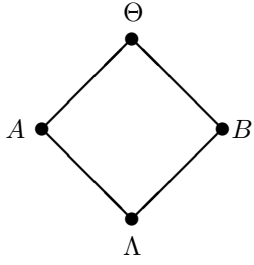
\begin{figure}
\centering
\setlength{\unitlength}{6mm}
\begin{picture}(6,6)\thicklines
\put(3,1){\circle*{0.3}}
\put(1,3){\circle*{0.3}}
\put(5,3){\circle*{0.3}}
\put(3,5){\circle*{0.3}}
\put(1,3){\line(1,1){2}}
\put(1,3){\line(1,-1){2}}
\put(5,3){\line(-1,1){2}}
\put(5,3){\line(-1,-1){2}}
\put(0.2,2.8){$A$}
\put(5.2,2.8){$B$}
\put(2.8,5.4){$\Theta$}
\put(2.8,0.2){$\Lambda$}
\end{picture}
\caption{Lattice Structure of $\A$\label{Fig.4.1}}
\end{figure}

Fig. \ref{Fig.4.1}, which is called the Hasse diagram, shows the lattice structure of $\A$.

\vskip 5mm

\vskip 5mm

Next, we consider ${\cal M}_{\mu}$. Define
$$
{\cal M}_{\mu}^k:={\cal M}_{k\mu_y\times k\mu_x}.
$$
Then we have a partition as
\begin{align}\label{4.5}
{\cal M}_{\mu}=\bigcup_{k=1}^{\infty}{\cal M}_{\mu}^k.
\end{align}
Now if $s=rt$, then
\begin{align}\label{4.6}
{\cal M}_{\mu}^r\otimes I_t\subset {\cal M}_{\mu}^s.
\end{align}
Then it is reasonable to pose an order on ${\cal M}_{\mu}$ as
\begin{align}\label{4.7}
{\cal M}_{\mu}^r\prec {\cal M}_{\mu}^s,\quad \mbox{if}~ r|s.
\end{align}
Then we have the lattice structure on ${\cal M}_{\mu}$ as follows.

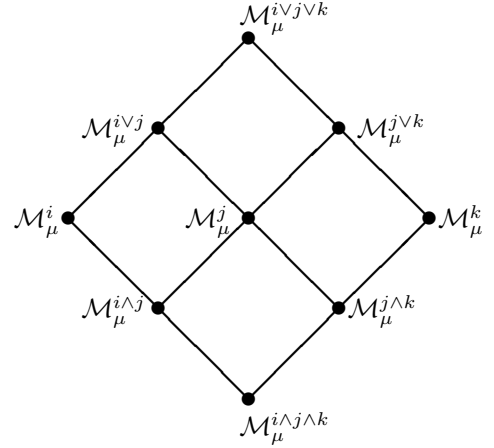
\begin{figure}
\centering
\setlength{\unitlength}{6mm}
\begin{picture}(10,10)(-0.5,0)\thicklines
\put(5,1){\circle*{0.3}}
\put(3,3){\circle*{0.3}}
\put(7,3){\circle*{0.3}}
\put(1,5){\circle*{0.3}}
\put(5,5){\circle*{0.3}}
\put(9,5){\circle*{0.3}}
\put(3,7){\circle*{0.3}}
\put(7,7){\circle*{0.3}}
\put(5,9){\circle*{0.3}}
\put(5,1){\line(-1,1){4}}
\put(5,1){\line(1,1){4}}
\put(3,3){\line(1,1){4}}
\put(7,3){\line(-1,1){4}}
\put(5,9){\line(-1,-1){4}}
\put(5,9){\line(1,-1){4}}
\put(4.8,0.2){${\cal M}_{\mu}^{i\wedge j\wedge k}$}
\put(1.3,2.8){${\cal M}_{\mu}^{i\wedge j}$}
\put(7.2,2.8){${\cal M}_{\mu}^{j\wedge k}$}
\put(-0.2,4.8){${\cal M}_{\mu}^{i}$}
\put(3.6,4.8){${\cal M}_{\mu}^{j}$}
\put(9.2,4.8){${\cal M}_{\mu}^{k}$}
\put(1.3,6.8){${\cal M}_{\mu}^{i\vee j}$}
\put(7.4,6.8){${\cal M}_{\mu}^{j\vee k}$}
\put(4.8,9.3){${\cal M}_{\mu}^{i\vee j\vee k}$}
\end{picture}
\caption{Lattice Structure of ${\cal M}_{\mu}$\label{Fig.4.2}}
\end{figure}

\begin{prp} \label{p4.7} ${\cal M}_{\mu}$ with the partial order determined by (\ref{4.7}) is a lattice, where
$$
\begin{array}{l}
\inf\left({\cal M}_{\mu}^p,{\cal M}_{\mu}^q\right)={\cal M}_{\mu}^{p\wedge q},\\
\sup\left({\cal M}_{\mu}^p,{\cal M}_{\mu}^q\right)={\cal M}_{\mu}^{p\vee q}.\\
\end{array}
$$
\end{prp}
Fig. \ref{Fig.4.2} shows the lattice structure of ${\cal M}_{\mu}$.

\vskip 5mm

\vskip 5mm

\section{Quotient Spaces of ${\cal M}_{\mu}$}

For statement ease, this section focuses on $\Gamma=\Gamma^1$, $\sim=\sim^{\Gamma^1}_{\ell}$, and $\Sigma_{\mu}=\Sigma^{\ell}_{\mu}$. We leave the parallel discussion for other quotient spaces to reader.

\subsection{Topology on Quotient Space}

We consider the topology of the quotient space  $\Sigma_{\mu}={\cal M}_{\mu}/\sim$. First, the topology of ${\cal M}_{\mu}$ should be clarified.

\begin{dfn}\label{d5.1} A natural topology on ${\cal M}_{\mu}$, denoted by ${\bf N}$, is constructed as follows.
\begin{itemize}
\item[(i)] For each ${\cal M}_{\mu}^k$, the natural topology of Euclidean space $\R^{k^2\mu_y\mu_x}$ is posed on ${\cal M}_{\mu}^k$, $k=1,2,\cdots$.
\item[(ii)] Each ${\cal M}_{\mu}^k$, $k=1,2,\cdots$, is considered as a clopen set of ${\cal M}_{\mu}$.
\end{itemize}
\end{dfn}

Then it is clear that under this natural topology the ${\cal M}_{\mu}$ has following property. We refer to any classical text book for the related concepts, e.g., \cite{dug66}.

\begin{prp}\label{p5.2} The topological space  $\left({\cal M}_{\mu},{\bf N}\right)$ is a Hausdorff space. It is also second countable.
\end{prp}

Next, we consider the quotient spaces of ${\cal M}_{\mu}$ with respect to equivalences.

\begin{dfn}\label{d5.3} Some quotient spaces of ${\cal M}_{\mu}$ under equivalences are defined as follows:
\begin{align}\label{5.1}
\begin{array}{l}
\Sigma_{\mu}^{\ell}:={\cal M}_{\mu}/\sim^{\varGamma^1}_{\ell};\\
\Sigma_{\mu}^{r}:={\cal M}_{\mu}/\sim^{\varGamma^1}_{r};\\
\Theta_{\mu}^{\ell}:={\cal M}_{\mu}/\sim^{\varGamma^2}_{\ell};\\
\Theta_{\mu}^{r}:={\cal M}_{\mu}/\sim^{\varGamma^2}_{r};\\
\end{array}
\end{align}
\end{dfn}

We pose the quotient topology on quotient spaces. To be specific, we focus on $\Sigma_{\mu}:=\Sigma_{\mu}^{\ell}$ and denote by $\left<A\right>:=\left<A\right>_{\ell}$.

Consider $\Sigma_{\mu}$. Assume $\pi:{\cal M}_{\mu}\ra \Sigma_{\mu}$ is the natural projection, i.e., $\pi:A\mapsto \left<A\right>$. Using the natural topology ${\bf N}$ on ${\cal M}_{\mu}$, then the quotient topology on $\Sigma_{\mu}$ deduced by the mapping $\pi$ is denoted by ${\bf Q}$. Precisely speaking, the open setc on $\Sigma_{\mu}$ is
$$
{\bf Q}:=\{\pi^{-1}(O)\;|\;O\in {\bf N}\}.
$$

Next, we consider the fiber bundle structure with ${\cal M}_{\mu}$ and $\Sigma_{\mu}$.

\begin{dfn}\label{d5.4} \cite{hus94} Let $T$ and $B$ be two topological spaces. $\pi: T\ra B$ is a continuous surjective mapping. Then
\begin{align}\label{5.2}
T\xrightarrow{~~~\pi~~~} B
\end{align}
is called a fiber bundle, with $T$ as the total space $B$ as the base space. For each $b\in B$, $\pi^{-1}(b)$ is called the fiber over $b$.
\end{dfn}

Then we have the following fiber bundle structure over ${\cal M}_{\mu}$ and $\Sigma_{\mu}$.

\begin{prp}\label{p5.5} Consider $\left({\cal M}_{\mu}, {\bf N}\right)$ as total space and $\left(\Sigma_{\mu},{\bf Q}\right)$ as base space, and 
$$
\pi:A\mapsto \A,\quad A\in {\cal M}_{\mu}
$$
is the natural projection. Then
\begin{align}\label{5.3}
\left({\cal M}_{\mu}, {\bf N}\right)\xrightarrow{~~~\pi~~~}\left(\Sigma_{\mu},{\bf Q}\right)
\end{align}
is a fiber bundle. The fiber of $\A$ is $\{B\;|\; B\in \A\}$.
\end{prp}

\subsection{Vector Space Structure on Quotient Space}

\begin{dfn}\label{d5.6} Consider $\Sigma_{\mu}$.
\begin{itemize}
\item[(i)] Let $\left<A\right>,\left<B\right>\in \Sigma_{\mu}$.
\begin{align}\label{5.4}
\begin{array}{l}
\left<A\right>\lplus \left<B\right>:= \left<A\lplus B\right>,\\
\left<A\right>\lminus \left<B\right>:= \left<A\lminus B\right>.
\end{array}
\end{align}
\item[(ii)] Let $\left<A\right>\in \Sigma_{\mu}$, $r\in \R$.
\begin{align}\label{5.5}
r\cdot\left<A\right>:= \left<rA\right>.
\end{align}
\end{itemize}
\end{dfn}

The following proposition can easily be verified by a straightforward computation.

\begin{prp}\label{p5.7} The addition $\lplus$ and scalar product $\cdot$ defined by (\ref{5.4}) and (\ref{5.5}) for  $\Sigma_{\mu}$ are properly defined. Moreover, $\left(\Sigma_{\mu},\lplus,\cdot\right)$ is a vector space.
\end{prp}

\begin{dfn}\label{d5.8} Consider $\Sigma_1$. Assume $\A,\B\in \Sigma_1$, then
\begin{align}\label{5.6}
\A \ltimes \B:=\left<A\ltimes B\right>.
\end{align}
\end{dfn}

It is easy to verify that (\ref{5.6}) is properly defined. That is, the right hand side of (\ref{5.6}) is independent of the choice of representatives $A\in \A$ and $B\in \B$.

\begin{dfn}\label{d5.9} \cite{gre81} A vector space $V$ (over $\R$) with a binary operator (called a product) $*:V\times V\ra V$
is called an algebra, if the following distributive law is satisfied.
$$
\begin{array}{l}
x*(ay_1+by_2)=ax*y_1+bx*y_2,\\
(ay_1+by_2)*x=ay_1*x+by_2*x.
\end{array}
$$
\end{dfn}

Then the following proposition is obvious.

\begin{prp}\label{p5.9} $\Sigma_1$ with the product define by  (\ref{5.6}) is an algebra.
\end{prp}

The followings are some fundamental properties of the algebra $(\Sigma_1,\ltimes)$.

\begin{dfn}\label{d5.8} Consider $\Sigma={\cal M}/\sim_{\ell}$.
\begin{itemize}
\item[(i)]
\begin{align}\label{5.6}%5.301}
\A^T:=\left<A^T\right>,\quad \A\in \Sigma.
\end{align}
\item[(ii)] Let $A\in {\cal M}_1$. If $A$ is invertible (where $\Gamma=\Gamma^1$), then
\begin{align}\label{5.7}% 302}
\A^{-1}:=\left<A^{-1}\right>,\quad \A\in \Sigma.
\end{align}
\item[(iii)] Let $\A,\B\in \Sigma_1$. $\A$ and $\B$ are said to be equivalent, if (where $\Gamma=\Gamma^1$) there exist nonsingular $P$ and $Q$, such that
\begin{align}\label{5.8}%5.303}
\left<P\right>\ltimes \A \left<Q\right> =\B.
\end{align}
\item[(iv)] Let $\A,\B\in \Sigma_1$. $\A$ and $\B$ are said to be congruent,  if (where $\Gamma=\Gamma^1$) there exists nonsingular $P$, such that
\begin{align}\label{5.9}%5.304}
\left<P\right>^T\ltimes \A \ltimes\left<P\right> =\B.
\end{align}
\item[(v)] Let $\A,\B\in \Sigma_1$. $\A$ and $\B$ are said to be similar,  if (where $\Gamma=\Gamma^1$) there exists a nonsingular $P$, such that
\begin{align}\label{5.10}%5.304}
\left<P\right>^{-1}\ltimes \A \ltimes\left<P\right> =\B.
\end{align}
\end{itemize}
\end{dfn}

\begin{rem}\label{r5.9} Using Definition \ref{d5.8}, many related properties of matrices can be extended naturally to $\Sigma_{\mu}$. For example:
\begin{itemize}
\item[(i)]
\begin{align}\label{5.11}
\left(\A \ltimes \B\right)^T=\B^T\ltimes \A^T.
\end{align}
\item[(ii)]If $A,B\in {\cal M}_1$ are invertible, then
\begin{align}\label{5.12}
\left(\A \ltimes \B\right)^{-1}=\B^{-1}\ltimes \A^{-1}.
\end{align}
\item[(iii)] If $\A$ is equivalent (congruent, similar) to $\B$, and  $\B$ is equivalent (correspondingly, congruent, similar) to $\left<C\right>$, then $\A$ is equivalent (correspondingly, congruent, similar) to $\left<C\right>$.
\end{itemize}\end{rem}

\subsection{Analytic Functions on $\Sigma_1$}

\begin{prp}\label{p5.10} Let $\A\in \Sigma_1$. Then
\begin{align}\label{5.13}
e^{\A}=\left<e^A\right>.
\end{align}
\end{prp}

\begin{IEEEproof} We have only to prove that if $A\sim B$, then $e^A\sim e^B$. Let $\inf(A,B)=\Theta$. It is enough to prove that $e^A\sim e^{\Theta}$. Note that there exists an identity matrix, say $I_k$, such that $A=\Theta\otimes I_k$. Then it is enough to show that
\begin{align}\label{5.14}
e^{\Theta\otimes I_k}=e^{\Theta}\otimes I_k.
\end{align}
Assume $\Theta\in {\cal M}_{n\times n}$, then we have
$$
W(\Theta\otimes I_k)W^{-1}=I_k\otimes \Theta=\diag\underbrace{(A,A,\cdots,A)}_k,
$$
where $W=W_{[n,k]}$.\footnote{$$W_{[m,n]}=\left(I_n\otimes \d_m^1, I_n\otimes \d_m^2,\cdots,I_n\otimes \d_m^m\right)$$ is called the swap matrix.

Let $A\in {\cal M}_{m\times n}$ and $B\in {\cal M}_{p\times q}$. Then \cite{che12}
$$
W_{[m,p]}(A\otimes B)W_{[q,n]}=B\otimes A.
$$
}
Then we have
$$
\begin{array}{ccl}
e^{A}&=&e^{W^{-1}(I_k\otimes \Theta)W}\\
~&=&W^{-1}e^{\diag(\Theta,\Theta,\cdots,\Theta)}W\\
~&=&W^{-1}\diag(e^{\Theta},e^{\Theta},\cdots,e^{\Theta})W\\
~&=&W^{-1}\left(I_k\otimes e^{\Theta}\right)W,\\
~&=&e^{\Theta}\otimes I_k.
\end{array}
$$
\end{IEEEproof}

\begin{rem}\label{r5.11} Eq. (\ref{5.13}) is of particular importance. In the sequel, it will be used to construct dimension-free general linear group (DFGLG) $\GL(\R)$ from dimension-free general linear algebra (DFGLA) $\gl(\R)$.
\end{rem}

In general, an analytic function can be defined over $\Sigma_1$:

\begin{prp}\label{p5.12} Let $f$ be an analytic function. Then $f(\A)=\left<f(A)\right>$. 
\end{prp}

\begin{IEEEproof} Note that
$$
\A\ltimes \B=\left<A\ltimes B\right>,
$$
and
$$
\A\lplus \B=\left<A\lplus B\right>.
$$
Using Taylor expansion, we have
$$
\begin{array}{ccl}
f(\A)&=&\lplus_{i=0}^{\infty}\frac{f^{(n)}(0)}{n!}\A^n\\
~&=&\lplus_{i=0}^{\infty}\frac{f^{(n)}(0)}{n!}\left<A^n\right>\\
~&=&\left<\lplus_{i=0}^{\infty}\frac{f^{(n)}(0)}{n!}A^n\right>\\
~&=&\left<f(A)\right>.
\end{array}
$$
\end{IEEEproof}

\section{Lie Algebraic Structure on $\Sigma_1$}

As aforementioned, $\left(\Sigma_1,\lplus\right)$ is a vector space, the Lie bracket on $\Sigma_1$ can be defined as follows:
\begin{align}\label{6.2}
[\left<A\right>,\left<B\right>]:=\left<A\right>\ltimes\left<B\right>\lminus
\left<B\right>\ltimes\left<A\right>,\quad \A,\B\in \Sigma_1.
\end{align}

It is easy to verify the following lemma.

\begin{lem}\label{l6.1} Let $A\in {\cal M}_1^r$, $B\in {\cal M}_1^s$, $C\in {\cal M}_1^t$, and $n=r\vee s\vee t$. Then
\begin{align}\label{6.3}
\begin{array}{ccl}
(A\ltimes B)\ltimes C&=&A\ltimes (B\ltimes C)\\
~&=&\left(A\otimes I_{n/r}\right)\left(B\otimes I_{n/s}\right)\left(C\otimes I_{n/t}\right).
\end{array}
\end{align}
\end{lem}

Using Lemma \ref{l6.1}, one can prove the following result by a straightforward computation.

\begin{prp}\label{p6.2} $\Sigma_1$ with Lie bracket defined by (\ref{6.2}) is a Lie algebra,
called the DFGLA, denoted by $\gl(\R)$ (or $\gl(\C)$ as the vector space is considered over $\C$).
\end{prp}

In the following, we consider only $\gl(\R)$, the results are all available for $\gl(\C)$, unless elsewhere is stated.

Define a function on ${\cal M}_{1}$ as follows.

\begin{dfn}\label{d6.3} Let $A\in {\cal M}_1^k$. Then
\begin{align}\label{6.4}
\Tr(A):=\frac{1}{k}\trace(A),
\end{align}
which is called the dimension-free trace.
\end{dfn}

\begin{rem}\label{r6.4}
\begin{itemize}
\item[(i)] It is easy to verify that if $A\sim^{\Gamma^1} B$ ($\sim^{\Gamma^1}$ stands for either
$\sim^{\Gamma^1}_{\ell}$ or $\sim^{\Gamma^1}_r$), then
\begin{align}\label{6.5}
\Tr(A)=\Tr(B).
\end{align}
Hence we can define
\begin{align}\label{6.6}
\Tr(\A):=\Tr(A),\quad \A\in \Sigma_1.
\end{align}
\item[(ii)] Eq. (\ref{6.4}) can be replaced by
\begin{align}\label{6.601}
\Tr(A):=\frac{1} {\trace(\Gamma_k)}\trace(A).
\end{align}
Then we have that if $A\sim^{\Gamma} B$, then
(\ref{6.5}) holds. Hence we can also define
\begin{align}\label{6.7}
\Tr(\A_{\Gamma}):=\Tr(A),\quad \A_{\Gamma}\in {\cal M}_{\mu}/\sim^{\Gamma}.
\end{align}
\end{itemize}
\end{rem}

\begin{dfn}\label{d6.201}
\begin{itemize}
\item[(i)] Let $(L,[\cdot,\cdot])$ be a Lie algebra, $V\subset L$ is a subspace of $L$. If 
$(V,[\cdot,\cdot])$ is also a Lie algebra, then it is called a Lie sub-algebra  of $L$. 
\item[(ii)] Assume $V\subset L$ is a Lie sub-algebra. $V$ is called an ideal of $L$, if
$$
[X, A]\in V, \quad X\in L,\;A\in V.
$$
\end{itemize}
\end{dfn}

\begin{dfn}\label{d6.202}
\begin{itemize}
\item[(i)] Let $(L,[\cdot,\cdot]_L)$ and $(H,[\cdot,\cdot]_H)$ be two Lie algebras. $\pi:L\ra H$ is called a homomorphism, if
\begin{align} \label{6.701}    
\pi([X,Y]_{L})=[\pi(X),\pi(Y)]_{H},\quad X,Y\in L.
\end{align}
\item[(ii)] Assume  $\pi:L\ra H$ is a Lie algebra homomorphism, if $\pi$ is a bijective mapping, $\pi$ is called an isomorphism.
\end{itemize}
\end{dfn}
    
In the following example some Lie sub-algebras of $\gl(\R)$ are constructed.
\begin{exa}\label{e6.3}
\begin{enumerate}
\item Let
\begin{align}\label{6.8}
{\cal A}:=\{\A\in \gl(\R)\;|\;Tr(\A)=0\}.
\end{align}
Then ${\cal A}$ is an ideal of $\gl(\R)$.
\item Given an $\left<M\right>\in \gl(\R)$, define
\begin{align}\label{6.9}
\begin{array}{l}
\gl(\left<M\right>,\R):=\{\A\in \gl(\R)\;|\;\\
~~\A \ltimes \left<M\right>\lplus \left<M\right>\ltimes \A^T=0\}.
\end{array}
\end{align}
Then $\gl(\left<M\right>,\R)$ is a Lie sub-algebra of $\gl(\R)$.
\begin{itemize}
\item[(i)] Assume
$$
M=\begin{bmatrix}
0&1\\1&0
\end{bmatrix},
$$
then the $\gl(\left<M\right>,\R)$ is called the orthogonal algebra.
\item[(ii)] Assume
$$
M=\begin{bmatrix}
0&1\\-1&0
\end{bmatrix},
$$
then the $\gl(\left<M\right>,\R)$ is called the symplectic algebra.
\item[(iii)] If $\left<M\right>$ and $\left<N\right>$ are congruence, then
$\gl(\left<M\right>,\R)$ and $\gl(\left<N\right>,\R)$ are isomorphic.
\end{itemize}
\end{enumerate}
\end{exa}

\section{Killing Form on $\gl(\R)$}

Next, we define the Killing form on $\Sigma_1=\gl(\R)$.

\begin{dfn}\label{d7.1}%d6.4}
\cite{hum72} Let $G$ be a lie algebra and $A\in G$. Define
$$
\ad_AX:=[A,X],\quad X\in G.
$$
Then $\ad_A\in \End(G)$ is an endomorphism. (Equivalently, if $\dim(G)=n$, $\ad_A\in \gl(n,\R)$.)
\end{dfn}

\begin{prp}\label{p7.2}  Let $A\in \gl(n,\R)$. Then $\ad_A\in \gl(n^2,\R)$ is expressed as
\begin{align}\label{7.1}%6.13}
\ad_A=I_n\otimes A-A^T\otimes I_n.
\end{align}
\end{prp}

\begin{IEEEproof}
Since
$$
\ad_AX=[A,X]=AX-XA.
$$
Express it into standard linear mapping form, we have \cite{che07}
$$
V_c(\ad_AX)=(I_n\otimes A-A^T\otimes I_n)V_c(X),
$$
where $V_c(A)$ is the column stacking form of $A$.
\end{IEEEproof}

Now we can define
$$
\ad: \Sigma_1\ra \End(\Sigma_1)
$$ as follows:
\begin{align}\label{7.2}
\ad_{\A}\left<X\right>:=[\A,\left<X\right>]=\left<[A,X]\right>.
\end{align}
Using Proposition \ref{p7.2}, we have
$$
\begin{array}{ccl}
\ad_{\A}\left<X\right>&=&[\A\left<X\right>]\\
~&=&\left<(A\otimes I_n-A^T\otimes I_n)X\right>\\
~&=&[\left<\ad_A\right>,\left<X\right>]\\
~&=&\left<\ad_A\right> \left<X\right>.
\end{array}
$$
Hence
\begin{align}\label{7.3}
\ad_{\A}=\left<\ad_A\right>.
\end{align}
Hence $\ad_{\A}$ itself can be considered as an element in $\gl(\R)$.
From this perspective, we have the following result.

\begin{prp}\label{p7.3}
$\ad: \gl(\R)\ra \gl(\R)$ is a homomorphism.
\end{prp}

\begin{IEEEproof}
$$
\begin{array}{l}
\ad_{\left<[A,B]\right>}\left<X\right>\\
~~=\left<[A,B]\right>\left<X\right>\\
~~=[\A\ltimes\B\lminus \B\ltimes \A]\left<X\right>\\
~~=\left<[A,B]X\right>\\
~~=\left(\ad_{\A}\ad_{\B}-\ad_{\B}\ad_{\A}\right)\left<X\right>\\
~~=[\ad_{\A},\ad_{\B}]\left<X\right>.\\
\end{array}
$$
Hence, we have
$$
\ad_{\left<[A,B]\right>}=[\ad_{\A},\ad_{\B}].
$$
\end{IEEEproof}

\begin{dfn}\label{d7.4}%d6.6}
\begin{itemize}
\item[(i)] Let $A,B\in {\cal M}_1$. Then the Killing form of $A$ and $B$ is define by
\begin{align}\label{7.4}%6.11}
\left(A,B\right):=Tr(\ad_A\ltimes\ad_B).
\end{align}
\item[(ii)] Let $\A,\B\in \Sigma_1$. Then the Killing form of $\A$ and $\B$ is define by
\begin{align}\label{7.5}%6.12}
\left(\A,\B\right):=\left(A,B\right).
\end{align}
\end{itemize}
\end{dfn}

We have to prove the Killing form (\ref{7.5}) is properly defined.

\begin{prp}\label{p7.5}%6.9}
The Killing form (\ref{7.5}) of $\Sigma_1$ is properly defined.
\end{prp}
\begin{IEEEproof}
Denote by $A_1\in \A$ and $B_1\in \B$ the minimum elements in $\A$ and $\B$ respectively. It is enough to prove that for any $A\in \A$ and $B\in \B$
\begin{align}\label{7.6}%6.14}
\left(A,B\right)=\left(A_1,B_1\right).
\end{align}

Assume $A_1\in {\cal M}_1^r$, $B_1\in {\cal M}_1^s$, and $r\vee s=t$.

Using formula (\ref{7.1}) and noticing the facts that (a) $\trace(AB)=\trace(BA)$; (b) $\trace(A^T)=\trace(A)$; (c) $\trace(A\otimes B)=\trace(A)\trace(B)$, we calculate that
$$
\begin{array}{l}
\Tr(\ad_{A_1}\ltimes \ad_{B_1})\\
~~=\Tr\left[\left(\ad_{A_1}\otimes I_{t^2/r^2}\right) \left(\ad_{B_1}\otimes I_{t^2/s^2}\right)\right]\\
~~=\Tr\left\{\left[\left(I_r\otimes A_1-A_1^T\otimes I_r\right)\otimes I_{t^2/r^2}\right]\right.\\
~~\left.\left[\left(I_s\otimes B_1-B_1^T\otimes I_s\right)\otimes I_{t^2/s^2}\right]\right\}\\
~~=\Tr\left\{\left[I_t\otimes \left(A_1\otimes I_{t/r}\right)-
\left(A^T_1\otimes I_{t/r}\right)\otimes I_t\right]\right.\\
~~\left.\left[I_t\otimes \left(B_1\otimes I_{t/s}\right)-
\left(B^T_1\otimes I_{t/s}\right)\otimes I_t\right]\right\}\\
~~=\Tr\left[I_t\otimes (A_1\ltimes B_1)+(A_1^T\ltimes B_1^T)\otimes I_t\right.\\
~~-\left(A_1^T\otimes I_{t/r}\right)\otimes  \left(B_1\otimes I_{t/s}\right)\\
~~\left.-\left(B_1^T\otimes I_{t/s}\right)\otimes \left(A_1\otimes I_{t/r}\right) \right]\\
~~=\frac{2}{t}\trace(A_1\ltimes B_1)-\frac{2}{rs}\trace(A_1)\trace(B_1).
\end{array}
$$
Finally, assume $A=A_1\otimes I_p$ and $B=B_1\otimes I_q$, then we have
$$
\begin{array}{l}
\Tr(\ad_{A}\ltimes \ad_{B})\\
~~=\Tr\left\{\left[\left(\ad_{A_1}\otimes I_{t^2/r^2}\right)\otimes I_p\right]\right.\\
~~\left.\left[\left(\ad_{B_1}\otimes I_{t^2/s^2}\right)\otimes I_q\right]\right\}\\
~~=\frac{2}{t}\trace(A_1\ltimes B_1)-\frac{2}{rs}\trace(A_1)\trace(B_1).
\end{array}
$$
\end{IEEEproof}

The following properties of Killing form on $\gl(\R)$ follow the definition immediately.

\begin{prp}\label{p7.6}
The Killing form on $\gl(\R)$ has the following propositions:
\begin{itemize}
\item[(i)]
\begin{align}\label{7.7}
\left(\A,\B\right)=\left(\B,\A\right),\quad \A,\B\in \gl(\R).
\end{align}
\item[(ii)]
\begin{align}\label{7.8}
\begin{array}{l}
\left(a\A_1\lplus b\A_2,\B\right)=a\left(\A_1,\B\right)\lplus b\left(\A_2,\B\right),\\
\qquad a,b\in \R.
\end{array}
\end{align}
\item[(iii)]
\begin{align}\label{7.9}
\left(\ad_{\A} \left<X\right>,\left<Y\right>\right)+\left(\left<X\right>,\ad_{\A} \left<Y\right>\right)=0.
\end{align}
\end{itemize}
\end{prp}

\section{Colculsion}

In this paper the DFES is proposed. The vector space structure is posed. Based on the equivalence, the quotient topology is constructed. Then the Lie bracket is defined, which turns the DFES a Lie algebra, which is of infinite dimension. It is then shown that this algebra has certain properties of finite dimensional Lie algebra. Particularly, the Killing form is constructed for it. Its properties are also investigated.

\end{document}